\documentclass[12pt,a4paper]{article}

\usepackage[english]{babel}
\usepackage[active]{srcltx}
\usepackage{amsmath,amssymb,amsthm,amsfonts,amscd}

\newcommand{\rav}{\stackrel{\triangle}{=}}
\newcommand{\upartial}{\partial}

\newcommand{\leqref}[1]{\stackrel{(\ref{#1})}{\leq}}
\newcommand{\geqref}[1]{\stackrel{(\ref{#1})}{\geq}}
\newcommand{\rref}[1]{$(\ref{#1})$}
\newcommand{\bi}[1]{{\it{ \bf #1}}}
\newcommand{\td}[1]{\widetilde{#1}}
\newcommand{\doc}{{\em{Proof }}}
\newcommand{\bo}{
\hfill {$\Box$}}
\theoremstyle{plain}
\newtheorem{theorem}{Theorem}[section]
\newtheorem{lemma}{Lemma}
\newtheorem{conjecture}{Conjecture}[section]
\newtheorem{corollary}{Corollary}[section]
\newtheorem{proposition}{Proposition}[section]

\def\theequation{\thesection.\arabic{equation}}
\theoremstyle{remark}
\newtheorem{remark}{Remark}

\theoremstyle{definition}
\newtheorem{definition}{Definition}

\begin{document}



\title{
Necessity of limiting co-state arc\\ in Bolza-type infinite horizon problem}

\author{D.V. Khlopin\\
 khlopin@imm.uran.ru\\
Krasovskii Institute of Mathematics and Mechanics\\ Yekaterinburg, Russia
}

\maketitle

\begin{abstract}
 We investigate necessary conditions of optimality for the Bolza-type infinite horizon problem
   with  free right end.
  The optimality is understood in the sense of weakly uniformly overtaking optimal control.
    No previous knowledge in the asymptotic behaviour of
    trajectories or adjoint variables
    is necessary.
    Following Seierstad’s idea, we obtain the necessary boundary condition at infinity in the form of
  a transversality condition for the maximum principle.
    Those transversality conditions may be expressed in the integral form through an Aseev--Kryazhimskii-type formulae for co-state arcs.
 The connection between these  formulae and limiting gradients of payoff function at infinity is identified; several conditions under which it is possible to explicitly specify the co-state arc through those Aseev--Kryazhimskii-type formulae are found.
  For      infinite horizon problem  of  Bolza type, an example is given to clarify the use of  the Aseev--Kryazhimskii formula  as explicit expression of the co-state arc.
\end{abstract}


{\bf 49K15;  49J52;  91B62}

$$\ $$

   The first necessary conditions of optimality for infinite-horizon control problems were proved \cite{ppp} on the verge of 1950--60s by L.S. Pontryagin and his associates (for the problems with the right end fixed at infinity).
       Only later \cite{Halkin} was the Maximum Principle proved for a reasonably broad class of problems, and yet the transversality-type conditions at infinity were not provided. A significant number
\cite{Halkin,kam,mich,smirn,ssbook} of such conditions
was proposed.
Thus, the Maximum Principle for infinite horizon was not complete, and the set of extremals obtained through it was too broad;  see
\cite{car,Halkin,mich,shell},\cite[Sect. 6]{kr_as},\cite[Example
10.2]{norv}.

      The principal obstacle on the way to transversality conditions at infinity is the fact that it is necessary to find  the asymptotic conditions on the adjoint equation (i.e., on the linear system) that would be satisfied by at least one but not by all of its solutions. It was first done in \cite{aucl} for linear autonomous control system  through passing to a functional space that allowed to extend all necessary solutions to infinity in the unique way. For the co-state arc, there was proved a formula that supplemented the Maximum Principle, making it a complete system. In the papers  \cite{kr_as3,kr_as,kab,av,av_new}, more general formula (the Aseev--Kryazhimskii formula) was proved  for certain other  classes of nonlinear  control problems. It takes the form of an improper integral of a function, the summability of which on the whole half-line is provided by means of imposing the asymptotic conditions (similar to the dominating discount conditions) on the system
      (for details, refer to Subsect.3.1 of this paper or \cite[Sect.16]{kr_as},\cite{av_new}, \cite[Sect.6]{MyJDC}).

       Another way to decrease the number of solutions of such an incomplete system of relations was proposed by Seierstad \cite{norv}. He considered a set of shortened problems and the corresponding expressions of the Maximum Principle, in each of which he obtained a co-state arc as a component of the solution of the corresponding  system. Under sufficiently strong assumptions he made, the co-state arc, obtained as a pointwise limit, satisfied the maximum principle of original infinite horizon problem. In \cite{MyJDC}, in line with this idea, the author obtained the necessary
  condition at infinity for the general class of control problems with fixed  initial state and free terminal state at infinity; a co-state arc was also expressed in the improper integral form, similar  Aseev--Kryazhimskii formula.

   This paper is concerned with the infinite horizon problem  of  Bolza type
   with  free right end.
  The optimal control is understood in the sense of weakly uniformly overtaking optimal control \cite{HamKen}.
   It
is assumed a priori that a uniformly weakly
overtaking optimal control exists (for discussion of existence
refer to \cite{bald,baum,car,dm,zaslav,z1}).
No previous knowledge in the asymptotic behaviour of
    trajectories or adjoint variables
    is necessary.

    A necessary boundary condition at infinity is found in the form of a
   transversality condition for the maximum principle (see Theorem~\ref{1}).
 A co-state arc that would satisfy this condition is constructed in the form of the pointwise limit of co-state arcs for shortened problems with tailored penalty function.
    The integral form of the transversality conditions may be seen as the limit of Aseev--Kryazhimskii-type formulae for co-state arcs (see \rref{1p},\rref{0p}).
  The connection between these  formulae and limiting subdifferential of payoff function at infinity is clarified
  (see Corollary~\ref{resume}).
  Several conditions are found (see Subsect. \ref{formulae}), under which the constructed transversality conditions complement the Maximum Principle, making it a complete system of relations; in this case, the  co-state arc may be  explicitly expressed by formulae \rref{1s},\rref{0s},\rref{1ss}, including
  Aseev--Kryazhimskii formula \rref{1sss}.
  These conditions are weaker than those from \cite{MyJDC,kab}; in particular, we dispensed with the necessity to assume the payoff function to be finite at infinity in
 \cite[Corollary 14]{MyJDC},\cite[Theorem 4]{kab}.
  The example from Subsect.~\ref{ex} serves to clarify the utility of the Aseev--Kryazhimskii formula as an explicit expression of the co-state arc in
     infinite horizon problems  of  Bolza type.

   Although the proof of Theorem~\ref{1} follows \cite{MyJDC}, it is not based on previous papers \cite{MyJDC,MyArkh,my2}, being significantly strengthened and streamlined; in particular, we no longer need to embed the set of admissible controls into a compact space of generalized controls.
   However, as in \cite{MyJDC}, the dynamics of the function is assumed to be smooth  in~$x$; it is important for the method used to construct the penalty function (see Appendix~A).

  \section{Problem statement and definitions}

  We consider the time interval
  ${\mathbb{T}}\rav
  {\mathbb{R}}_{\geq 0}.$
The phase space of the
 control system is the
finite-dimensional Euclidean space ${\mathbb{X}}\rav{\mathbb{R}}^m$.

  Consider the following optimal control problem
\begin{subequations}
\begin{eqnarray}
    \textrm{Minimize } l(b)+\int_{0}^\infty f_0 (x,u,t) dt  \label{sys0}\\
    \textrm{subject to } \dot{x}=f(x,u,t),\quad t>0,\quad u\in U(t)  \label{sys}\\
    x(0)=b\in{\mathcal{C}}.    \label{sysK}
\end{eqnarray}
\end{subequations}
  Here,  functions $l$ and $f_0$ are scalar; $x$ is the  state variable taking values in  ${\mathbb{X}};$ and $u$ is the control parameter.

Suppose that  $U$ is a map  from $T$ to the set of all closed subsets of a finite-dimensional Euclidean space~${\mathbb{U}}$, and the
     graph of $U$ is a Borel set. As for the class of admissible controls, we consider the set of measurable functions $u(\cdot)$ bounded for any time compact
 such that $u(t)\in U(t)$ holds for a.a. $t\in{\mathbb{T}}$.
 Denote the set of admissible controls
   by~${\mathfrak{U}}$.

 We assume the following conditions hold:
\begin{itemize}
  \item  ${\mathcal{C}}$ is a closed subset of ${\mathbb{X}}$;
  \item   $l$ is taken to be locally Lipshitz continuous on $x$;
  \item  $f$ is Borel measurable in $(t,u)$ and continuously differentiable
in $x$;
  \item   for each admissible control $u$, the map  $(t,x)\mapsto f(x,u(t),t)$ satisfies the
 sublinear growth condition (see, for example, \cite[1.4.4]{tovst1});
  \item   $f_0$ is measurable in $(t,u)$, continuously differentiable
in $x$,  and lower semicontinuous in $u$;
  \item  $\frac{\partial f}{\partial x},\frac{\partial f_0}{\partial x}$ are measurable in $(t,u)$ and
 locally Lipshitz continuous on  $x$.
\end{itemize}

Now we can assign the solution of \rref{sys} to each admissible control~$u$ and initial condition $b\in{\mathbb{X}}$.
The solution is unique and it can be extended to the whole ${\mathbb{T}}$.
Let
us denote it by~$x(b,u;\cdot).$
The pair $(b,u)$ will be called an admissible control process
if~$u\in {\mathfrak{U}},$ $b\in{\mathcal{C}}.$

   For every $(b, u)\in{\mathbb{X}}\times{\mathfrak{U}}$, assign the number
\[
   J(b,u;T)\rav\int_0^T f_0\big(x(b,u;t),u(t),t\big)dt
\]
to each $T\in{\mathbb{T}}$.

Note the integral in \rref{sys0} may not exist in the general case; if it exists, it could tend to $+\infty$ and may even do so for all admissible controls.
In such a situation, we have to adopt the special notion of optimality.
 A number of such definitions of optimality for control problems on infinite horizon were proposed \cite{car1,car,stern,slovak}.
We will essentially understand the objective of~\rref{sys0} to be
 \[ \textrm{Minimize } l(b)+J(b,u;T) \textrm{   for very large  } T,\]
 which conforms to the notion of optimality based on ``overtaking,'' specifically, to the notion of ``weakly uniform overtaking'';
(see \cite{brock1973,HamKen,car1}).
\begin{definition}
 If an admissible process $(b^*,u^*)\in{\mathcal{C}}\times{\mathcal{U}}$ satisfies
\[ \liminf_{T\to\infty}
  \Big(
    l(b^*)+J(b^*,u^*;T)
    -
      \inf_{(b,u)\in{\mathcal{C}}\times{\mathcal{U}}}\big(l(b)+J(b,u;T)\big)
     \Big)
      \geq 0,\]
call it
    a uniformly overtaking optimal process.
\end{definition}

Assume we know in advance the test sequence of times at which  the result of this process can be interesting;
let these times form an unboundedly increasing sequence of times $\tau_n$.

\begin{definition}\label{275}
  If an admissible process $(b^*,u^*)\in{\mathcal{C}}\times{\mathcal{U}}$ satisfies
\[ \liminf_{n\to\infty}
  \bigg(
    l(b^*)+J(b^*,u^*;\tau_n)
    -
      \inf_{(b,u)\in{\mathcal{C}}\times{\mathcal{U}}}\big(l(b)+J(b,u;\tau_n)\big)
     \bigg)
      \geq 0,\]
call it
    a uniformly $\tau$-overtaking optimal process.
\end{definition}
There also exists a criterion of weakly uniformly overtaking  optimality
(see \cite{brock1973}), which we can now reformulate in the following way:
\begin{definition}
  If an admissible process $(b^*,u^*)$ is a uniformly $\tau$-overtaking optimal process for some
   unboundedly increasing sequence of times $\tau_n$,
call it
    a weakly uniformly overtaking  optimal process.
\end{definition}

   Note that we can also do the optimization by way of choosing the proper sequence $\tau$.
   For instance,  we could allow the players of infinite-horizon differential games to choose their own test sequence in addition to their control, and consider the corresponding Nash equilibrium; let us call it a test Nash  equilibrium.
   Any test Nash equilibrium is weakly   uniformly overtaking optimal.
    Moreover, there are examples of infinite-horizon differential games with  a unique uniformly overtaking optimal equilibrium and with a unique test Nash equilibrium such that for each player
    the uniformly overtaking optimal equilibrium is the worst of all weakly   uniformly overtaking optimal  Nash equilibria, and the
    test Nash equilibrium is the best of  weakly   uniformly overtaking optimal  Nash equilibria
    (See \cite[Example 1]{MyArkh}).

  Hereinafter assume there exists an admissible process $(b^*,u^*)$ that is uniformly $\tau$-overtaking optimal for some unbounded increasing sequence of numbers $\tau_n.$
  Denote by $x^*$ the solution generated by this process.

  Note that this assumption is equivalent to existence of a weakly uniformly overtaking solution for problem \rref{sys0}--\rref{sysK}. This paper is not concerned with when is it so; for the issues of existence of solutions of infinite-horizon problems that are optimal with respect to one or another criterion, refer to \cite{bald,baum,car1,car,dm,tan_rugh,zaslav,z1}.

Let us also fix a sequence $\tau=(\tau_n)_{n\in{\mathbb{N}}}$.
Slightly simplifying the notation when passing from the sequence
$\tau$ to its subsequence $\tau'$, we will
plainly write ``subsequence $\tau'\subset\tau\,$.''

\section{PMP relations}
\setcounter{equation}{0}

Let us recall a few definitions from vector analysis \cite{cl_new},\cite[Sect.4]{vinter}.

Let $g:{\mathbb{X}}\to{\mathbb{R}}\cup\{+\infty\}$ be a lower semicontinuous function.
 A vector $\zeta\in{\mathbb{X}}$
is said to be a proximal subgradient of $g$  at $b$ provided that there exist a
neighborhood $\Omega$ of $x$ and a number $\sigma\geq 0$ such that
\[g(\xi) \geq g(b) + \zeta(\xi-b)-\sigma||\xi-b||^2\qquad \forall\xi\in \Omega.\]
The set of proximal subgradients at $x$ (which may be empty, and which is not
necessarily closed, open, or bounded but which is convex) is denoted $\upartial_P g(b)$,
and is referred to as the proximal subdifferential. This set  is nonempty for all $x$ in
a dense subset of $\{b\,|\,g(b)<+\infty\}.$
 Following
 \cite[Theorem 4.6.2(a)]{vinter},
  denote the limiting subdifferential of $g$ at $b$ by $\upartial_L g(b)$; it
 consists of all
$\zeta$ in ${\mathbb{X}}^*$  such that
\[
\exists \textrm{ sequences of }  b_n\in{\mathbb{X}},\zeta_n\in \upartial_P g(b_n),b_n\to b,
\zeta_n\to\zeta.\]
 Following
 \cite[Theorem 4.6.2(b)]{vinter},
 denote the singular limiting (asymptotic limiting) subdifferential of $g$ at $b$ by $\upartial^0_{L} g(b)$; it consists of all
$\zeta$ in ${\mathbb{X}}^*$  such that
\[
\exists \textrm{ sequences of }  b_n\in{\mathbb{X}},s_n\in{\mathbb{T}},\zeta_n\in \upartial_P g(b_n),b_n\to b,s_n\downarrow 0,
s_n\zeta_n\to\zeta.\]

If $g$ is Lipshitz continuous near $b$, then $\upartial_{L} g(b)$ is nonempty, moreover $co\, \upartial_L g(b)=\upartial_{Clarke} g(b),\upartial^0_{L} g(b)=\{0\}$
(see \cite[Sect. 4]{vinter}).

We say $\zeta\in{\mathbb{X}}$ is a proximal normal to  ${\mathcal{C}}$ at $b\in{\mathcal{C}}$ if there exists $\sigma\geq 0$ such that
\[\zeta(\xi-b)\leq \sigma||\xi-b||^2\qquad \forall\xi\in {\mathcal{C}}.\]
The set of such $\zeta$ is denoted by  $N_P^{{\mathcal{C}}}(b)$
and referred to as the proximal normal cone \cite[Ch~1.3]{clarke}.
 Denote the limiting normal cone  to ${\mathcal{C}}$ at $b$, by $N_L^{{\mathcal{C}}}(b)$, it consists of all
$\zeta$ in ${\mathbb{X}}$  such that
\[
\exists \textrm{ sequences of }  b_n\in{\mathcal{C}},\zeta_n\in N_P^{{\mathcal{C}}}(b_n),b_n\to b,\zeta_n\to\zeta.
\]

   Let us now proceed to the relations of the Pontryagin Maximum Principle.

  Let the Hamilton--Pontryagin function
  $H:{\mathbb{X}}\times {U}\times{\mathbb{X}}\times{\mathbb{T}}\times {\mathbb{T}}\mapsto{\mathbb{R}}$
  be given by
   $H(x,u,\psi,\lambda,t)\rav\psi f\big(x,u,t\big)-\lambda
   f_0\big(x,u,t\big).$
 Let us introduce the relations and boundary condition:
\begin{subequations}
 \begin{eqnarray}
   \label{sys_x}
       \dot{x}(t)&=& f\big(x(t),u(t),t\big);\\
   \label{sys_psi}
       -\dot{\psi}(t)&=&\frac{\partial
       H}{\partial x}\big(x(t),u(t),\psi(t),\lambda,t\big);\\
   \label{maxH}
\!\!\!\sup_{v\in
U(t)}H\big(x(t),v,\psi(t),\lambda,t\big)&=&
        H\big(x(t),u(t),\psi(t),\lambda,t\big);\\
   \label{dob}
   ||\psi(0)||+\lambda&=&1.
   \end{eqnarray}
\end{subequations}
 It is easily seen that, for each $u\in{\mathfrak{U}}$ for each initial
  condition, system~\rref{sys_x}--\rref{sys_psi} has a local solution,
   and each solution of these relations can be extended to the
    whole~${\mathbb{T}}$.

The PMP relations written above are well known, and their necessity for infinite horizon was first proved in \cite{ppp,Halkin}.
Note, however, that these relations are incomplete; in the general case, they must be supplemented with an additional condition, a boundary condition at infinity (for details, refer to \cite{kr_as,car,Halkin,mich,nnn,norv,ssbook}).

   As it was mentioned in papers \cite{norv,nnn}, under certain additional assumptions,
   the condition that some nontrivial solution of PMP is a pointwise limit of solutions of this system where the adjoint variable is zero at a sufficiently large time can serve as a necessary condition of optimality
  for control problems with the free right end.
    This condition was refined into the following definition (see definition of vanishing solution in \cite{MyJDC} ):
\begin{definition}
    A nontrivial solution
   $(\lambda^*,\psi^*)$ of  \rref{sys_x}--\rref{maxH}
   associated with $(x^*,u^*)$ is called $\tau$-limiting (or just limiting)
   if $(x^*,\psi^*,\lambda^*)$ is a pointwise limit
      of
        solutions $(x_n,\psi_n,\lambda_n)$ of boundary value problems
\begin{subequations}
 \begin{eqnarray}
   \label{sys_x_k}
       \dot{x}(t)&=& f\big(x(t),u^*(t),t\big);\\
        \label{sys_psi_k}
       -\dot{\psi}(t)&=&\frac{\partial
       H}{\partial x}\big(x(t),u^*(t),\psi(t),\lambda,t\big);\\
   \label{sys_l_k}
       \dot{\lambda}(t)&=&0;\\
   \label{dob_k}
   \psi_n(\tau'_n)&=&0 
 \end{eqnarray}
 \end{subequations}
  for some test subsequence
                     $\tau'\subset\tau$.
\end{definition}
\begin{remark}
   Since the right-hand side of~\rref{sys_psi_k}--\rref{sys_l_k} is homogeneous by $(\psi,\lambda)$,   without loss of generality we can say  that
      $\lambda_n+||\psi_n(0)||=\lambda^*+||\psi^*(0)||.$ Or, if $\lambda^*>0$, then without loss of generality we can say that
      $\lambda_n=\lambda^*$.
\end{remark}
 Note that the case of strong optimality allows us to strengthen this definition significantly \cite{my2}.

    The principal result of this paper is the following:
\begin{theorem}
\label{1}
    Let the process $(b^*,u^*)$ be  a uniformly $\tau$-overtaking process for problem \rref{sys0}--\rref{sysK}.

    Then
    for  $(b^*,u^*)$ there exists a $\tau$-limiting solution $(\psi^*,\lambda^*)$ of  system \rref{sys_x}-\rref{dob} such that
\begin{equation}\label{400}
    \psi^*(0)\in \lambda^*\upartial_L l(b^*)+N_L^{{\mathcal{C}}}(b^*).
\end{equation}
\end{theorem}
The proof of Theorem \ref{1} in relocated to section \ref{doc}.
The sketch of the proof is as follows:
Consider a sequence of control problems on the intervals of the form $[0,\tau_n]$; functional \rref{sys0} is supplemented with the penalty in the form of an integral factor that depends only on $t$ and  $u$.  The penalty may be chosen such that its sufficiently small value would guarantee that the trajectory of the control system will not exit the given integral funnel of solutions generated by the control~$u^*$; in Appendix A, we demonstrate the consruction of such a penalty function. By the Ekeland principle these problems have optimal solutions, for which we write the corresponding PMP relations. By passing to subsequence we can find both the limit of these solutions and the PMP relations that are satisfied for it. The definition of optimality implies that the penalty is small, whence the fact that the limiting solution could only be generated by the control $u^*$, in particular, $(x^*,u^*,\psi^*,\lambda^*)$ must satisfy \rref{sys_x}-\rref{dob}.
To prove that $(\psi^*,\lambda^*)$ is $\tau$-limiting, from the right-hands end of the  solutions that are optimal on $[0,\tau_n]$, we launch in reverse time the solutions  $(x_n,\psi_n,\lambda_n)$ that are generated by the optimal control~$u^*$.
They are contained in the integral funnel generated by~$u^*$ around $(x^*,\psi^*,\lambda^*)$ due to the fact that our construction of the penalty allows us to make the width of the funnel at the initial time as small as we want; therefore, the constructed sequence converges to $(x^*(0),\psi^*(0),\lambda^*)$ at the initial time. Thus, by theorem of continuous dependence on initial conditions, we have the required convergence on every finite interval.
\begin{proposition}
\label{1262}
  Under assumptions of the theorem, we can also demand that, in addition,
\begin{equation}
  ||J(x_n(0),u^*;\tau_n)-J(b^*,u^*;\tau_n)||\to 0\textrm{ as } n\to\infty\label{1214}.
  \end{equation}
\end{proposition}
See the proof of Proposition~\ref{1262} in  section \ref{manydoc}.

 Straight from the definitions, we get
\begin{remark}
\label{11}
    Let the process $(b^*,u^*)$ be  uniformly overtaking for problem \rref{sys0}--\rref{sysK}.

    Then,  for each unbounded sequence of $\tau_n>0$ there exists a $\tau$-limiting solution $(\psi^*,\lambda^*)$ of  system
     \rref{sys_x}-\rref{dob}  satisfying \rref{400}.
\end{remark}
\begin{remark}
\label{12}
    Let the process $(b^*,u^*)$ be   weakly uniformly overtaking for \rref{sys0}--\rref{sysK}.

    Then,  for some unbounded sequence of  $\tau_n>0$  there exists a $\tau$-limiting solution $(\psi^*,\lambda^*)$ of  system
    \rref{sys_x}-\rref{dob}  satisfying \rref{400}.
\end{remark}

 \section{Formulae as corollaries of Theorem~\ref{1}}
\setcounter{equation}{0}
 \subsection{Limiting formulae.}

Let us use the fact that system
 \rref{sys_psi_k}--\rref{sys_l_k}  is linear. Denote by
${\mathbb{L}}$ the linear space of all real $m\times m$ matrices; here $m=dim\,{\mathbb{X}}$.
For each $\xi\in{\mathbb{X}}$, there exists a solution ${A}(\xi;t)\in C({\mathbb{T}}, {\mathbb{L}})$ of the Cauchy problem
 \begin{equation*}
 \frac{d{A}(\xi;t)}{dt} =\frac{\partial f }{\partial x}
 \big(x(\xi,u^*;t),u^*(t),t\big)
  A(\xi;t),\quad A(\xi;0)=1_{\mathbb{L}}.
\end{equation*}
Let us introduce the vector-valued function $I$ of time by the following rule: for every
  $T\in{\mathbb{T}}$,
  \begin{equation*}
  I(\xi;T)\rav\int_0^T
   \frac{\partial f_0}{\partial x}
    \big(x(\xi,u^*;t),u^*(t),t\big)
\, A(\xi;t)
  \,dt.
\end{equation*}
  Now, for each solution $(x,\psi,\lambda)$ of system \rref{sys_x_k}--\rref{sys_l_k}, for every $t\in{\mathbb{T}}$, we have the Cauchy formula
  \begin{equation}
   \label{4A}
   \psi(t)=\big(\psi(0)+\lambda I(x(0);t)\big)A^{-1}(x(0);t).
\end{equation}

  Note that if $\psi(\tau_n)=0$,
  then $\psi(t)=\lambda\big( I(x(0);t)-I(x(0);\tau_n)\big)A^{-1}(x(0);t);$
  in particular,
  $\psi(0)=-\lambda I(x(0);\tau_n).$
  Replacing $(x,\psi,\lambda)$ with the triples $(x_n,\psi_n,\lambda_n)$ from the definition of $\tau$-limiting solution, we obtain the following:
    \[\psi_n(0)=-\lambda_n I(x_n(0);\tau_n). \]
  By passing, if necessary, to the subsequence, we can guarantee that either the sequence of
  $I(x_n(0);\tau_n)$ converges or the sequence  of
  $I(x_n(0);\tau_n)/||I(x_n(0);\tau_n)||$ converges.
  Thus we prove
\begin{subequations}
\begin{proposition}
\label{affp} {
     Let the conditions of Theorem~\ref{1} be satisfied.
         Let there exist a $\tau$-limiting solution $(\lambda^*,\psi^*)$.

    Then, for some subsequence $\tau'\subset\tau$
    and a sequence of
    $\xi_n\in{\mathbb{X}}$ converging to $b^*$,
  accurately to a positive factor,
  one of the following two relations
also hold:
    \begin{equation}
    \label{1p}
\mbox{either}\quad  \lambda^*=1,\quad\psi^*(0)=  -\lim_{n\to\infty}
I(\xi_n;\tau'_n)\in \upartial_L l(b^*)+N_L^{{\mathcal{C}}}(b^*),
\end{equation}
      and the sequence of $I(\xi_n;\tau'_n)$ has the finite limit,
    \begin{equation}
    \label{0p}
  \mbox{\ \ \ \ or}\quad  \lambda^*=0,\quad  {\psi^*(0)}=
  -\lim_{n\to\infty}\frac{I(\xi_n;\tau'_n)}{||I(\xi_n;\tau'_n)||}\in N_L^{{\mathcal{C}}}(b^*),
  \end{equation}
   this limit is valid, the sequence of $||I(\xi_n;\tau'_n)||$ is unbounded.
  }
\end{proposition}
\end{subequations}

For the case $I(\xi;\cdot)\equiv I(b^*;\cdot)$,
we have the explicit formulae
\begin{subequations}
\begin{corollary}
\label{affs} {
Let, under assumptions of Theorem~\ref{1}, $f$ and $f_0$ be linear by $x$.
Let there exist a $\tau$-limiting solution $(\lambda^*,\psi^*)$ corresponding to $(b^*,x^*)$.

    Then, for some subsequence $\tau'\subset\tau$,
  accurately to a positive factor,
  one of the two following relations
also holds:
    \begin{equation}
    \label{1s}
 \mbox{either}\quad  \lambda^*=1,\quad \psi^*(0)= -\lim_{n\to\infty} I(b^*;\tau'_n)\in \upartial_L l(b^*)+N_L^{{\mathcal{C}}}(b^*),
  \end{equation}
   and the sequence of $I(b^*;\tau'_n)$ has the finite limit
    \begin{equation}
    \label{0s}
    \mbox{\ \ \ \ or}\quad  \lambda^*=0,\quad {\psi^*(0)}=
  -\lim_{n\to\infty}\frac{I(b^*;\tau'_n)}{||I(b^*;\tau'_n)||}\in N_L^{{\mathcal{C}}}(b^*),
  \end{equation}
   this limit is valid, the sequence of $||I(b^*;\tau'_n)||$ is unbounded.
   }
\end{corollary}
\end{subequations}

 In some circumstances we can prove more. We can prove that the formulae above point out the unique $\tau$-limiting solution, for example,
\begin{corollary}
\label{affss} {
  Let the conditions of Theorem~\ref{1} be satisfied.
         Let $(\lambda^*,\psi^*)$ be some $\tau$-limiting solution corresponding to $(b^*,x^*)$.
Let there be
\begin{equation}
   \label{I1} \displaystyle I_*=\lim_{n\to\infty,\xi\to b^*}  I(\xi;\tau_n)\in{\mathbb{X}}
\end{equation}

\begin{subequations}
    Then,
  accurately to a positive factor, there exists a unique
    $\tau$-limiting solution $(\lambda^*,\psi^*)$;
    this solution is given by the following rules:
    \begin{eqnarray}
    \lambda^*=1,\qquad \psi^*(0)=-I_*\in \upartial_L l(b^*)+N_L^{{\mathcal{C}}}(b^*),\nonumber\\
    \psi^*(T)=
     \Big(-I_*+\int_{0}^T
   \frac{\partial f_0}{\partial x}\big(x^*(t),u^*(t),t\big)A(b^*;t)
  \,dt\Big)\,A^{-1}(b^*,T)\quad\forall T\geq 0.
    \label{1ss}
    \end{eqnarray}

    Moreover, if $I_*$ from \rref{I1} is independent of~$\tau$
         (i.e. $\displaystyle I_*=\lim_{t\to\infty,\xi\to b^*} I(\xi;t)$ holds),  we obtain
           the Aseev--Kryazhimskii formula
    \begin{eqnarray}
    \label{1sss}
    \psi^*(T)=
     -\int_{T}^\infty
   \frac{\partial f_0}{\partial x}\big(x^*(t),u^*(t),t\big)\, A(b^*;t)
  \,dt\,A^{-1}(b^*,T)\qquad\forall T\geq 0,
  \end{eqnarray}
 where the improper integral is understood in the Riemann sense. In particular,
\[     -\int_{0}^\infty
   \frac{\partial f_0}{\partial x}\big(x^*(t),u^*(t),t\big)\, A(b^*;t)
  \,dt\in\upartial_L l(b^*)+N_L^{{\mathcal{C}}}(b^*)\]
\end{subequations}
}
\end{corollary}
Note that formula~\rref{1sss}, as an explicit expression for the co-state arc in control problems on infinite horizon with the free right end, is  sufficiently studied. For some classes of linear problems, it was already obtained in paper~\cite{aucl}.
For the case of a singleton ${\mathcal{C}}$,
it was extended to a part of stationary problems by S.M.~Aseev and A.V.~Kryazhimskii in papers \cite{kr_as3,kr_as,kab};
for a number of nonstationary control systems, S.M.~Aseev and V.M.~Veliov proved its necessity in papers
   \cite{av,av_new}. Note that
   in  papers \cite{av,av_new},  \rref{1sss} is a necessary condition of local weakly overtaking optimality.
   It can be shown that for the problem considered, the necessary condition from \cite[Theorem~3.1]{sagara},\cite[Theorem 2.1]{norv} also reduces to~\rref{1sss}.
Let us also note that under assumptions of those papers, the improper integral in \rref{1sss} necessarily converges in the Lebesgue sense.
In addition to the papers cited above, the conditions that are sufficient to use~\rref{1sss} were studied by the author in \cite{MyJDC}. For details, see the references in \cite[Sect.~6.3]{MyJDC},\cite[Sect.~13]{kr_as}.

In the case of a strongly optimal solution \cite{car1}, we can guarantee the feasibility of one of the formulas \rref{1s},\rref{0s} as a necessary condition of optimality (for a singleton~${\mathcal{C}}$, refer to \cite{my2}).
This is not true in case of a uniformly overtaking process, the example of which will be stated below.
Let us also note that, as exhibited by the example offered in \cite[Ex.2]{MyJDC} (a modification of the Halkin example \cite{Halkin}),
  in Theorem \ref{1} we can not
replace the $\tau$-overtaking optimality (weakly uniformly overtaking optimality, uniformly overtaking optimality) with the DH-optimality (weekly agreeable, agreeable optimality;
\cite{car1}).

   Note that formula \rref{1ss} uses the information about the subsequence of $\tau.$ The usefulness of this fact in the search for Nash equilibria is showcased by the game from~\cite{MyArkh}. In this seemingly antagonistic linear game, a Nash equilibrium is chosen by players choosing the subsequences $\tau$.

Denote ${\mathcal{T}}\rav\{\tau_n\,|\,n\in{\mathbb{N}}\}.$
For a differentiable function $g:{\mathbb{X}}\times{\mathbb{T}}\to{\mathbb{X}}$,
similarly to the definitions of limiting subdifferential and singular limiting subdifferential,
let us introduce the generalized subdifferential of $g$ at infinity along $\tau$ by the following rule:
\begin{eqnarray*}\upartial^{1}_L g(b,\infty_\tau)&\rav&
  \{\zeta\,|\,\exists \textrm{ sequences of }  b_n\in{\mathbb{X}},t_n\in{\mathcal{T}},
  \zeta_n\in \upartial_P g(b,t_n),\\ &\ &b_n\to b,
 t_n\to\infty,\zeta_n\to\zeta\}.
\end{eqnarray*}
  Since in the general case it may be empty, let us also introduce a
  singular subdifferential in the following way:
\begin{eqnarray*}\upartial^{0}_L g(b,\infty_\tau)&\rav&
  \{\zeta\,|\,\exists \textrm{ sequences of }  b_n\in{\mathbb{X}},t_n\in{\mathcal{T}},s_n\in{\mathbb{T}},\zeta_n\in \upartial_P g(b,t_n),\\
  &\ &b_n\to b,t_n\to\infty,s_n\to 0,
s_n\zeta_n\to\zeta\}.
\end{eqnarray*}

Note that for all $T\in{\mathbb{T}}$,
   \[\frac{\partial J}{\partial x}(x,u^*;T)\equiv I(x;T),\qquad \upartial_P (-J)(x,u^*;T)\subset\upartial_L (-J)(x,u^*;T)\equiv\{-I(x;T)\}.\]
   Moreover,
   for each $T\in{\mathbb{T}}$,
     $\upartial_P (-J)(x,u^*;T)=\{-I(x;T)\}$
 for     a dense subset of ${\mathbb{X}}.$

   Then, the statements above could be summed up as
\begin{corollary}
\label{resume}
Let the conditions of Theorem~\ref{1} be satisfied.
         Then there exists a $\tau$-limiting solution  $(\lambda^*,\psi^*)$ corresponding to $(b^*,x^*)$.

Provided $\upartial^{0}_L J(b^*,u^*;\infty_\tau)=\{0\}$, we have
\begin{description}
  \item [a)] $\lambda^*>0$,
  \item [b)] $\psi^*(0)\in \lambda^*\upartial_L l(b^*)+N_L^{{\mathcal{C}}}(b^*),$
  \item [c)] $\psi^*(0)\in \lambda^*\upartial^{1}_L (-J)(b^*,u^*;\infty_\tau).$
  \item [d)] $\upartial^{1}_L (-J)(b^*,u^*;\infty_\tau)=\{-I_*\}$ iff \rref{I1} holds.
            \end{description}

Provided $\upartial^{1}_L J(b^*,u^*;\infty_\tau)=\varnothing$, we have
\begin{description}
  \item [e)] $\lambda^*=0$,
  \item [f)] $\psi^*(0)\in N_L^{{\mathcal{C}}}(b^*),$
  \item [g)] $\psi^*(0)\in \upartial^{0}_L (-J)(b^*,u^*;\infty_\tau)\setminus\{0\}.$
           \end{description}

Passing from~$\tau$ to its subsequence, we can always provide one of these statements.
\end{corollary}

  Consider a $p$ in $s\upartial^{1}_L (-J)(b^*,u^*;\infty_\tau)$ for some $s>0$.
 Then, for some subsequence $\tau'\subset\tau$ and a sequence of $b_n$ converging  to $b^*$,
the sequence of $p_n\in s\upartial_P (-J)(b_n,u^*;\tau'_n)=\{-I(b_n;\tau'_n)\}$
  converges to $p$. Thus, $-sI(b_n;\tau'_n)\to p$ as $n\to\infty.$
  Take a solution $(x_n,\psi_n,s)$ of \rref{sys_x_k}--\rref{sys_l_k} with
  $x_n(0)=b_n,\psi_n(0)=-sI(b_n;\tau'_n),\lambda=s,$
  and a solution $\psi^*$ of \rref{sys_psi} with
  $x\equiv x^*,u=u^*,\lambda=s.$
  By  the Cauchy formula \rref{4A}, we have $\psi_n(\tau'_n)=0,$ i.e. \rref{dob_k}.
   By  theorem on continuous dependence of the solutions of a differential equation,
   the solutions $(x_n,\psi_n,\lambda_n)$  converges pointwise to the
 solution $(x^*,\psi^*,s)$. Then,
 the solution
   $(s,\psi^*)$ is $\tau$-limiting  if $(s,\psi^*)$ is solution of \rref{sys_x}--\rref{maxH}.
  Thus, the nontrivial solution
   $(s,\psi^*)$ of  \rref{sys_x}--\rref{maxH}
   associated with $(x^*,u^*)$ is limiting if $\psi^*(0)\in s\upartial^{1}_L (-J)(b^*,u^*;\infty_\tau), s>0.$

   Similarly, nontrivial solution
   $(0,\psi^*)$ of  \rref{sys_x}--\rref{maxH}
   associated with $(x^*,u^*)$ is limiting if $\psi^*(0)\in\upartial^{0}_L (-J)(b^*,u^*;\infty_\tau).$

   We prove
\begin{theorem}
\label{2}
{A solution} $(\psi^*,\lambda^*)$ {of} \rref{sys_x}--\rref{maxH} {associated with} $(x^*,u^*)$
      {is  $\tau$-limiting iff}
      {a solution} $(\psi^*,\lambda^*)$  {of}  \rref{sys_x}--\rref{maxH} {associated with} $(x^*,u^*)$
      satisfies
\begin{eqnarray*}
{\textrm{either}}\quad
 \lambda^*>0,&\ &\psi^*(0)\in \lambda^*\upartial^{1}_L (-J)(x^*(0);u^*,\infty_{\tau});\\
{\textrm{or}}\quad
 \lambda^*=0,&\ &\psi^*(0)\in \upartial^{0}_L (-J)(x^*(0);u^*,\infty_\tau)\setminus\{0\}.
\end{eqnarray*}
\end{theorem}

\begin{corollary}
\label{3}
    Let the process $(b^*,u^*)$ be  a uniformly $\tau$-overtaking process for problem \rref{sys0}--\rref{sysK}.

    Then,
    for  $(b^*,u^*)$, there exists a nontrivial
     solution $(\psi^*,\lambda^*)$ of  system \rref{sys_x}-\rref{dob} such that
\begin{eqnarray*}
    &\ &\lambda^*\in\{0,1\};\\
    &\ &\psi^*(0)\in \upartial^{\lambda^*}_L (-J)(x^*(0);u^*,\infty_\tau);\\
    &\ &\psi^*(0)\in \lambda^*\upartial_L l\left(x^*(0)\right)+N_L^{{\mathcal{C}}}\left(x^*(0)\right).
\end{eqnarray*}
\end{corollary}

 In case of $\lambda>0$, form \rref{1ss} fits snugly with
    the interpretation of
  $\psi^*(0)$ as a marginal cost ("shadow price") associated with perturbing the initial condition. For infinite interval, such interpretation was discussed in the framework of the Ramsey model
in T.C. Koopmans \cite{coopmans}; under sufficiently broad assumptions, it was rigorously proved in \cite{aucl}.
For the latest results, refer to \cite{sagara}.

The mapping $b\mapsto J(b,u^*;t)$ is differentiable for all $t>0$,
thus, $\upartial_L J(b,u^*;t)$ is always a singleton set. Consequently, by passing form $\tau$ to a subset ${\mathcal{T}}'\subset{\mathcal{T}}$, we can always provide for exactly one of the sets
\begin{eqnarray*}
  \{\zeta\,|\,\exists \textrm{ sequences of }  t_n\in{\mathcal{T}}',
  \zeta_n\in \upartial_L J(b,u^*;t_n),
 t_n\to\infty,\zeta_n\to\zeta\}\\
  \{\zeta\,|\,\exists \textrm{ sequences of }  t_n\in{\mathcal{T}}',\zeta_n\in \upartial_L J(b,u^*;t_n),  t_n\to\infty,||\zeta_n||\to\infty,
\frac{\zeta_n}{||\zeta_n||}\to\zeta\}
\end{eqnarray*}
 to be a singleton set, i.e., we can provide the existence of one of the limits of \rref{1s},\rref{0s}.
 However, in general, we can not replace $\upartial^{1}_L$, $\upartial^{0}_L$ (respectively) in {Corollary}~\ref{resume} with them, see the example.

\subsection{Example}
\label{ex}

It is possible to construct an instance \cite[Example~3]{MyJDC} of the problem with fixed initial state, where
there is no nontrivial Lagrange multiplier $(1,\psi^*)$ (or $(0,\psi^*)$)
associated with $(x^*,u^*)$  that satisfies
  formula~\rref{1s} (or \rref{0s}).
  In this example for each limiting solution, i.e., for each limit of solutions of
  system~\rref{sys_x_k}-\rref{sys_l_k}, $x_n\not\equiv x^*$ for all $n\in{\mathbb{N}}.$
  Let us construct a Bolza type control problem by modifying this example.

Define a continuous monotonically nondecreasing convex scalar function~$f$ through the rule
\[
f(x)\rav
 \left\{
 \begin{array} {rcl}        0,       &\mathstrut&     x< 0,\\
                            \frac{x^2}{2},
                                         &\mathstrut&    0\leq x\leq 1,\\
                               x-\frac{1}{2},          &\mathstrut&    x>1.
 \end{array}            \right.
 \]
 Consider the following optimal control problem
\begin{subequations}
\begin{eqnarray}
    \textrm{Minimize } \int_{0}^\infty e^{-2t}x(x^4-5) dt   \label{ex1}\\
    \textrm{subject to } \dot{x}=f(x)+u,\quad t>0,\quad u\in[0,1], \label{ex2}\\
    x(0)\in{\mathcal{C}}\rav[-1,2].\label{ex3}
\end{eqnarray}
\end{subequations}

Consider an admissible process $b^*=0,u^*\equiv 0$,
then $x^*\equiv 0$.
It is possible to show (see Appendix~B) that in the present example
\begin{conjecture}
\label{924}
 The process $(b^*,u^*)=(0,0)$ is uniformly overtaking optimal   for problem \rref{ex1}-\rref{ex3}.
\end{conjecture}

 Note that, by $u^*\equiv 0$ and \rref{maxH}, for each solution $(\psi,\lambda)$ of PMP corresponding to $(b^*,u^*)$, we have $\psi(\cdot)\leq 0.$
 Now, by $x\equiv 0$, \rref{sys_psi} is equivalent to the equation
 \[\dot{\psi}=-5\lambda e^{-2t};\]
 thus, for each nontrivial solution $(\psi,\lambda)$ of the PMP, we have $\lambda>0.$
 In our case, $\lambda^*\upartial_L l(0)+N_L^{{\mathcal{C}}}(0)=\{0\}$ implies that, for any unbounded increasing sequence $\tau$, each $\tau$-limiting solution $(\psi^*,\lambda^*)$ satisfying Theorem~\ref{1} also satisfies $\psi^*(t)=\frac{5\lambda^*}{2}(e^{-2t}-1)<0$ for all $t>0.$
  Thus, a $\tau$-limiting solution of PMP satisfying Theorem~\ref{1} is unique, accurately to a positive factor.

  However, $A(0;\cdot)\equiv 1$,  $I(0;T)=-\int_{0}^T 5e^{-2t} dt<0$; then, there exists
\[I_{*}=\lim_{t\to\infty} I(0;T)<0.\]
Now, were formula~\rref{1s} to hold, $\psi$ would be positive, and by \rref{maxH} $u^*(t)>0$ for sufficient small $t>0,$ which is inconsistent with~$u^*\equiv 0$.
   Therefore, formula \rref{1s} does not apply here; thus,
  in this example,
there is no limiting  solution of the PMP such that its sequence of solutions $(x_n,\psi_n,\lambda_n)$ of \rref{sys_x_k}-\rref{sys_l_k} satisfies $x_n\equiv x^*.$

  And so, even if
  \begin{itemize}
                   \item  process $(b^*,u^*)$ is uniformly overtaking optimal;
                   \item  any nontrivial solution of the PMP that corresponds to it has $\lambda>0$
  (the problem is normal);
                   \item  a $\tau$-limiting solution of the PMP satisfying Theorem~\ref{1} is unique, accurately to a positive factor, and independent of the choice of $\tau$;
                   \item  the improper integral from \rref{1sss} exists in the Lebesgue sense;
                   \item  $b^* \in int\, {\mathcal{C}};$ $f_0$ is smooth,
                 \end{itemize}
 the formula $\psi^*(0)=-\lambda^* I_*$ (i.e., \rref{1s}) may not give a solution of PMP, in particular, may not give a $\tau$-limiting solution.

 This means that, for perturbed initial conditions, it is insufficient to compare the obtained gain $J(b,u^*;t)$ with the ``calculated'' gain $J(b^*,u^*;t)$ as $b\to b^*,t\to\infty$.
 In the general case, it is necessary, for $b$ sufficiently close to $b^*$, to move $J(\xi,u^*;t)$ when $\xi$ is near $b$, thereby finding $J(b,u^*;t)$, and afterwards make the perturbation $b-b^*$ as $t\to\infty$; this limiting set would contain~$\psi^*(0)$.
 One could say that $-\psi^*(0)$ must be  the limiting subgradient of $J(b,u^*;t)$ at $b^*$ along $\tau.$

\subsection{Formulae applicability.}
\label{formulae}

 The formulae \rref{1sss},\rref{1ss} (or even \rref{0s}) are the necessary conditions of optimality if
 either $\upartial^{1}_L J^*(b^*,\infty_\tau)$ or $\upartial^{0}_L J^*(b^*,\infty_\tau)\setminus\{0\}$ is a singleton set. In this case we
 have exactly one limiting PMP solution that corresponds to the optimal process.
  Let us obtain the more elegant (with respect to \rref{I1}) conditions that guarantee it.

\begin{proposition}
\label{s3_dd_kab1}
For each sequence of points $b_n$ that converges to~$b^*$, let
  \begin{equation}
   \label{620}
   I(b_n;\tau_n)-I(b^*;\tau_n)\to 0 \textrm{ as }
   J(b_n,u^*;\tau_n)-J(b^*,u^*;\tau_n)\to 0.
\end{equation}
Then, the results of Corollary~\ref{affs} hold.
 \end{proposition}

It is sometimes possible to verify condition~\rref{620} with ease. For example,

\begin{proposition}
\label{s3_dd_kab800}
 Assume that there exists  an open neighborhood $G\subset{\mathcal{C}}$ of $b^*$
 such that  the maps $b\mapsto I(b,\tau_n)-I(b^*,\tau_n)$ are equicontinuous on $G$
 for all $n\in{\mathbb{N}}.$
Suppose that $l$ is continuously differentiable in $x$ at $b.$

Then, the results of Corollary~\ref{affss} hold for $\lambda^*=1,I_*=-\frac{\partial l}{\partial x}(b^*)$.

 In addition, if the process $(b^*,u^*)$ is uniformly overtaking optimal
  for problem \rref{sys0}--\rref{sysK},  we obtain \rref{1sss}; in particular, there exists a unique solution $(\psi^*,1)$ of  system
     \rref{sys_x}-\rref{dob} satisfying $\psi^*(0)=-I_*$;
     this solution is $\tau$-limiting for any
    unbounded sequence of positive numbers $\tau_n$, and
 \begin{equation}
 \label{900}
    \int_{0}^{T}
   \frac{\partial f_0}{\partial x}\big(x^*(t),u^*(t),t\big)A(b^*;t)\,dt\to -\frac{\partial l}{\partial x}(b^*)\quad\mathrm{ as }\quad
   T\to\infty.
 \end{equation}
 \end{proposition}
Note that, as demonstrated by the example above, we can not omit the uniform continuity of $I(b,\tau_n)-I(b^*,\tau_n)$ from the assumptions of the Proposition.

  The results \cite[Corollary 14]{MyJDC},\cite[Theorem 4]{kab} show the necessity of
  \rref{1ss} if, for a fixed $b=b^*$, at large~$t$ the variance of $I$ is uniformly bounded from above by the value of variance of~$J$ for $u^*$ that are sufficiently close to the control $u$.
 For the Bolza functional, it is more natural to have such an estimate for fixed control, and not for a position, i.e., a function of differences $I(b;t)-I(b^*;t),$ $J(b,u^*;t)-J(b^*,u^*;t);$ let us do it.
\begin{proposition}
\label{s3_dd_kab}
Consider a function
$\omega\in C({\mathbb{T}}\times{\mathbb{T}},{\mathbb{T}})$ with $\omega(0,0)=0$.
For points $b$ from a neighborhood of $x^*$, let
  \begin{equation}
   \label{2742}
   ||I(b;\tau_n)-I(b;\tau_k)||
\leq
  \omega(1/\tau_k,|J(b,u^*;\tau_n)-J(b,u^*;\tau_k)|\big)
\end{equation}
 for all $n,k\in{\mathbb{N}} (n>k)$.

Then, the pair $(x^*,u^*)$ is normal,   and there exists
a limit $\displaystyle
 I_*=\lim_{n\to\infty}I(b^*;\tau_n)\in{\mathbb{X}}$;
 moreover, there exists
 a $\tau$-limiting solution $(1,\psi^*)$ such that
 \rref{1ss} and $\psi^*(0)=I_*$ hold.
 \end{proposition}

The proofs of all Propositions were placed in Section~\ref{manydoc}.


\section{Proof of Theorem~\ref{1}
 }
\label{doc}
\setcounter{equation}{0}

\subsection{Step 1:
  {\bi Choosing the metric}~$\rho$.}

  Set
   $E\rav{\mathbb{X}}\times{\mathbb{X}}\times {\mathbb{R}}.$
  Let ${\mathcal{S}}$ be a ball in $E$ centered at $y_0\rav(b^*,0,0)$ with the radius $2.$
Let the mapping $a:E\times {\mathbb{U}}\times {\mathbb{T}}$ be the right-hand side of the system
\begin{subequations}
 \begin{eqnarray}
   \label{sys_x_}
       \dot{x}(t)&=& f\big(x(t),u,t\big);\\
        \label{sys_psi_}
       \dot{\psi}(t)&=&-\frac{\partial
       H}{\partial x}\big(x(t),u,\psi(t),\lambda,t\big);\\
   \label{sys_l_}
       \dot{\lambda}(t)&=&0.
\end{eqnarray}
\end{subequations}
 Note that
  the system \rref{sys_x_k}--\rref{sys_l_k} satisfies all the requirements we demand from a system \rref{a}.
  Let us now construct, similarly to Appendix~A, mappings $w,\rho$ for such system~$a$
  with  designated control $u^*$ and compact set ${\mathcal{S}}$.

  Thus,   ${\mathcal{C}}\times{\mathfrak{U}}_{\tau_n}$ is metrizable by  $\left((b,u),(\xi,v)\right)\mapsto ||b-\xi|| + \rho(u,v,\tau_n)$
 for each $n\in{\mathbb{N}}$, here
\[{\mathfrak{U}}_{\tau_n}\rav \{u\in{\mathfrak{U}}\,|\,u(t)=u^*(t).\]

\subsection{Step 2:
  {\bi Constructing the auxiliary optimal solution sequence}.}

It is easy to prove that
  $(b^*,u^*)$ is a uniformly $\tau$-overtaking optimal process iff
  there exists  a sequence of positive numbers $\gamma_n\in(0,1/2]$ converging to
zero such that
 \begin{equation}
 \label{1588}
     l(b)+J(b,u;\tau_n)+\gamma^2_n>l(b^*)+J(b^*,u^*;\tau_n)
  \qquad \forall n\in{\mathbb{N}}, (b,u)\in{\mathcal{C}}\times{\mathcal{U}}.
 \end{equation}
 In particular, for each $n$, the mapping
  $(b,u)\mapsto l(b)+J(b,u;\tau_n)$
  is bounded from below by ${\mathcal{C}}\times{\mathfrak{U}}_{\tau_n}.$

  Now, the Ekeland principle \cite[Theorem~5.3.1,(i)]{ekeland},
  \cite[Theorem~2.1.3]{conv_new},\cite[Theorem~7.5.1]{clarke} furnishes a process $(b_n,u_n)\in{\mathcal{C}}\times{\mathfrak{U}}_{\tau_n}$
  that minimizes
  \begin{subequations}
\begin{eqnarray}
    l(b)+{J}(b,u;\tau_n)+\gamma_n \big(||b-b_n|| + \rho(u,u_n,\tau_n)\big)
     \label{sys01}\\
    \textrm{subject to } \dot{x}=f(x,u,t),\quad t>0,\quad u\in U(t),  \label{sys1}\\
    x(0)=b\in{\mathcal{C}}.    \label{sysK1}
   \end{eqnarray}
\end{subequations}
  Moreover (see \cite[Theorem~5.3.1,(i)]{ekeland}),\cite[Theorem~2.1.3,(ii)]{conv_new}),
  \begin{subequations}
\begin{eqnarray}\label{to_w_}
l(b_n)+{J}(b_n,u_n;\tau_n)+
\gamma_n \big(||b^*-b_n|| + \rho(u^*,u_n,\tau_n)\big)
\leq
l(b^*)+{J}(b^*,u^*;\tau_n).\ \
\end{eqnarray}
Hence, by virtue of \rref{1588}, we have
\begin{eqnarray}
 ||b^*-b_n||<\gamma_n\to0\textrm{ as }n\to\infty\label{1048},&\ &\\
  \rho(u^*,u_n,\tau_n)<\gamma_n\to 0\textrm{ as }n\to\infty\label{1050}.&\ &
   \end{eqnarray}
\end{subequations}
 Define $\td{x}_n$ by the rule $\td{x}_n(t)\rav x(b_n,u_n;t)$ for all $t\geq 0.$
 Now, from \rref{1048}, we have  $\td{x}_n(0)=b_n\to b^*=x^*(0).$
Note that, by \rref{1050} and Lemma~A.1, $u_n$ converges in measure to $u^*$ on the whole ${\mathbb{T}}.$
Passing to the subsequence if necessary, we can say that
$u_n$ converges to $u^*$ a.e. on ${\mathbb{T}}.$

\subsection{Step 3:
  {\bi PMP for auxiliary solutions.}}

  Since $u_n$ provides a minimum of problem \rref{sys01}-\rref{sysK1}, it can, if need be, yield the Pontryagin Maximum Principle~\cite[Theorem~5.1.1]{cl_new}.

  Let    the function ${H}_n:{\mathbb{X}}\times{U}\times{\mathbb{X}}\times{\mathbb{T}}\times {\mathbb{T}}\mapsto{\mathbb{R}}$
 be given by
  \[\displaystyle {H}_n(x,u,\psi,\lambda,t)\rav \psi f(x,u,t) -\lambda f_0(x,u,t)-
  \lambda\gamma_n w(u^*(t),u,t).\]
  Then, by the  Maximum
Principle,
 there exist
 ${\lambda}_n\in(0,1]$, $\td{\psi}_n\in C({\mathbb{T}},{\mathbb{X}})$
such that
relation~\rref{dob} and
 the transversality conditions
 \begin{subequations}
\begin{eqnarray}
\td\psi_n(\tau_n)&=&0,\label{trans_n_max_}\\
\td\psi_n(0)&\in& \lambda_n\upartial_L l(b_n)+\lambda_n\gamma_n\zeta+N_L^{{\mathcal{C}}}(b_n)
\textrm{ for some }\zeta\in{\mathbb{X}} (||\zeta||\leq 1)\label{trans_0_max_}
   \end{eqnarray}
 hold,  and
 \begin{eqnarray}
   \label{sys_max_}
   \!\!\sup_{v\in
        U(t)}{H}_n\big(\td{x}_n(t),v,\td{\psi}_n(t),{\lambda}_n,t\big)&=&
           {H}_n\big(\td{x}_n(t),u_n(t),\td{\psi}_n(t),{\lambda}_n,t\big),\\
-\dot{\psi}_n(t)&=&\frac{\partial
{H}_n}{\partial
x}\big(\td{x}_n(t),u_n(t),\td{\psi}_n(t),{\lambda}_n,t\big)
\label{1575}
 \end{eqnarray}
\end{subequations}
   also hold for a.e. $t\in[0,\tau_n]$.
  Since $H_n-H$ is independent of  $x$, we can now see that \rref{1575} implies
    \rref{sys_psi_} for $(\td{\psi}_n,\lambda_n)$.

\subsection{Step 4:
  {\bi PMP for overtaking optimal solution}}
  Consider $\td{y}_n\equiv (\td{x}_n,\td{\psi}_n,{\lambda}_n)$ for each $n\in{\mathbb{N}};$ note that this is a solution of system \rref{sys_x_}--\rref{sys_l_}.

  Passing, if need be, to a subsequence, we can consider the subsequence of
  ${\lambda}_n \in (0,1]$ to converge to a certain $\lambda^*\in[0,1]$ and a subsequence of $\td{\psi}_n(0)$ to converge to some ${\psi}^*_0\in{\mathbb{X}}$ as well.
 By continuity, from \rref{dob}, we obtain the relation $||{\psi}^*_0||+{\lambda}^*=1.$
  We now have, for sufficiently large~$n$,
\begin{equation}\label{1208}
 ||\td{y}_n(0)-y_0||\leq \lambda^*+||\psi^*_0||+||\td{x}_n(0)-b^*||\leqref{dob} 1+||\td{x}_n(0)-b^*||<2,
\end{equation}
   i.e.,
  $\td{y}_n(0)\to(b^*,{\psi}^*_0,\lambda^*)\in int\,{\mathcal{S}}.$

   In addition, for each $T>0$ starting with some $n$, we have $\rho(u^*,u,T)\leq\rho(u^*,u,\tau_n)$ for all $u\in{\mathfrak{U}}$; now, from  \rref{1050}, we have $\rho(u^*,u_n,T)\to 0.$ Therefore, by Lemma~A.3, in every compact set, the subsequence of $\td{y}_n$ uniformly converges to the solution $y^*$
  of system \rref{sys_x_}--\rref{sys_l_} generated by the control $u^*$, i.e., to the solution of
  \rref{sys_x_k}--\rref{sys_l_k}. Moreover, $y^*(0)=(b^*,\td{\psi}^*_0,\lambda^*).$
  But then there exists a solution ${\psi}^*$ of equation \rref{sys_psi_k} such that
   $y^*=({x}^*,{\psi}^*,{\lambda}^*)$ with ${\psi}^*(0)={\psi}^*_0.$

 Passing to the limit in \rref{sys_max_}, we have condition \rref{maxH} for  $({x}^*,{\psi}^*,{\lambda}^*)$ for almost every $t>0$.
 Thus, the limit $({x}^*,{\psi}^*,{\lambda}^*)$ satisfies system \rref{sys_x}--\rref{dob} for $u=u^*$, i.e., system \rref{sys_x_k}--\rref{sys_l_k}.

Since $x_n(0)\to x^*(0)=b^*$, $\psi_n(0)\to\psi^*(0)$, $\gamma_n\to 0$ as $n\to\infty$, by passing to the limit in \rref{trans_0_max_} we also gain \rref{400}, i.e.,
\[\psi^*(0)\in \lambda^*\upartial_L l (b^*)+N_L^{{\mathcal{C}}}(b^*).\]
 We found the solution $({x}^*,{\psi}^*,{\lambda}^*)$ of all relations of the PMP. It remains to prove that it is $\tau$-limiting.

\subsection{Step 5:
  {\bi Backtracking}}

Since \rref{1050} and \rref{1208} imply that
$\rho(u^*,u_n,\tau_n)<\gamma_n<1/2<dist(\td{y}_n(0),bd\,{\mathcal{S}}),$
 and
$\td{y}_n\to y^*$, $\gamma_n\to 0$ as $n\to\infty$, we know that Lemma~A.2 guarantees
\begin{equation}\label{1227}
   \varkappa(\td{y}_n(\tau_n),\tau_n)\to y^*(0).
 \end{equation}

From the position $\td{y}_n(\tau_n)$, launch in reverse time a trajectory $y_n$ of system
\rref{sys_x_}--\rref{sys_l_} with the help of the control $u^*$.
 Then, $y_n(0)=\varkappa(\td{y}_n(\tau_n),\tau_n)$ (see \rref{1667}).
   Note that $y_n=(x_n,\psi_n,\lambda_n)$ satisfies
   \rref{sys_x_k}--\rref{sys_l_k}, and
   $\psi_n(\tau_n)=\td{\psi_n}(\tau_n)=0$.
It remains to prove that the solution $y^*=(x^*,\lambda^*,\psi^*)$ of all PMP relations is the limit (in the compact-open topology) of the trajectories
    $y_n=(x_n,\psi_n,\lambda_n).$
    The same fact also follows from \rref{1227} by the theorem on continuous dependence of the solution of a differential equation on initial conditions.

The proof of Theorem~\ref{1} is complete.

\begin{remark}
 By $x_n(\tau_n)=\td{x}_n(\tau_n)$, it is safe to say that
 $\left(\tau_n,x_n(\tau_n)\right)$ belongs to the reachability domain of system \rref{sys} from the domain ${\mathcal{C}}.$
\end{remark}
\begin{remark}
 On the other hand, we may not generally say that $x_n(0)\in{\mathcal{C}},$
  for ${\mathcal{C}}=\{0\}$, the counterexample was shown in \cite{MyJDC}.
\end{remark}
\begin{remark}
    The result of this theorem remains valid if we omit the condition of Lipshitz continuity on $x$ for the functions
    $f_0,\frac{\partial f}{\partial x},\frac{\partial f_0}{\partial x}.$
    It is sufficient to assume the mapping $(t,x)\mapsto\frac{\partial f}{\partial x}\left(t,x,u^*(t)\right)$
    to be locally Lipshitz continuous in a domain~$G$ that contains the trajectories generated by $u^*$ (solutions of \rref{sys_x_k}) that are sufficiently close to $x^*$.
\end{remark}

\section{Proofs of propositions}
\label{manydoc}
\doc {\it of {Proposition}}~\ref{1262}.
\setcounter{equation}{0}

  Let us describe the necessary amendments to the proof of Theorem~1.
  As a first step, we need to broaden system \rref{sys_x_}--\rref{sys_l_}.
 Let us supplement it with the equation
\begin{equation}
  \label{sys_g_}
       \dot{\sigma}(t)=f_0\big(x(t),u,t\big).
 \end{equation}
  Now,
   $E\rav{\mathbb{X}}\times{\mathbb{X}}\times {\mathbb{R}}\times {\mathbb{R}},$
   ${\mathcal{S}}$ is a ball in $E$ centered at $(b^*,0,0,0)$ with the radius $2.$
  Let us now define $w,\rho$ for system \rref{sys_x_}--\rref{sys_l_},\rref{sys_g_}.

  On step~4, we introduce  $\td{\sigma}_n$ for each $n\in{\mathbb{N}}$ by the rule
  $\td{\sigma}_n(t)=J(\td{x}_n(0),u_n;t)$ $(\forall t\geq 0).$
  Now, $\td{y}_n\rav(\td{x}_n,\td{\psi}_n,\td{\lambda}_n,\td{\sigma}_n).$
  Estimate \rref{1208} remains valid by $J(\td{x}_n(0),u_n;0)=0$.
  By $\td{y}_n(0)\to (b^*,\psi^*_0,\lambda^*,0)$ and Lemma~A.3,
 the sequence of $\td{y}_n$ converges to the solution $y^*=({x}^*,{\psi}^*,{\lambda}^*,{\sigma}^*),$
 and ${\sigma}^*(t)=J(b^*,u^*;t)$ for all $t>0.$

  On step~5, let us likewise launch in the reverse time a trajectory $y_n$ of new system
\rref{sys_x_}--\rref{sys_l_},\rref{sys_g_} with the help of the control $u^*$ from the position $\td{y}_n(\tau_n)$.
  Then, along with $x_n,\psi_n,\lambda_n$, we will get some $\sigma_n.$
  It is easy to check that
  \[\sigma_n(t)-\td{\sigma}_n(t)=
  J(b_n,u_n;\tau_n)-J(x_n(0),u^*;\tau_n)\qquad \forall t\geq 0.\]
  Besides, by $\td\sigma_n(0)=J(b_n,u_n;0)=0$, \rref{1227} yields
  \[|J(b_n,u_n;\tau_n)-J(x_n(0),u^*;\tau_n)|=|\sigma_n(0)|\leq \gamma_n\to 0.\]

  By optimality of $u^*,u_n$ for their respective problems (see \rref{1588},\rref{to_w_}), we also have
  \[|J(b_n,u_n;\tau_n)-J(b^*,u^*;\tau_n)|\leq |l(b_n)-l(b^*)|+\gamma_n.\]
  Adding the last two estimates, by the triangle inequality, we obtain
  \rref{1214}.
  \bo

    \doc {\it of Proposition}~\ref{s3_dd_kab1}.

By Theorem \ref{1}, for this problem there exists the
$\tau$-limiting solution
  $(\lambda^*,\psi^*)$.
   In particular, there exists a sequence of solutions $(x_n,\psi_n,\lambda_n)$  of system \rref{sys_x_}--\rref{sys_l_} that converges to   $(x^*,\psi^*,\lambda^*)$ and satisfies \rref{dob_k}
  for some test subsequence
                     $\tau'\subset\tau$.
Moreover, by Proposition \ref{1262}, we can consider a subsequence of $J(x_n(0),u^*;\tau_n)-J(b^*,u^*;\tau_n)$ to tend to zero as well.
 Since $x_n(0)\to b^*$, the assumption of the corollary implies now
 $||I(x_n(0);\tau_n)-I(b^*;\tau_n)||\to 0.$

  Let us set $\lambda=\lambda_n, x=x^*$ in equation \rref{sys_psi_k} and consider its solution $\bar\psi_n$ such that
 $\bar\psi_n(\tau'_n)=0$.
 Now, $(x^*,\bar\psi_n,\lambda_n)$ satisfies \rref{sys_x_k}--\rref{sys_l_k},\rref{dob_k}.
 By \rref{4A} and \rref{620}
 \begin{equation}\label{631}
 \bar\psi_n(0)-\psi_n(0)=\lambda_n \big(I(x_n(0);\tau'_n)- I(x^*(0);\tau'_n)\big)\to 0.
 \end{equation}
 Moreover, $\psi_n\to \psi^*$ now implies $\bar\psi_n(0)\to \psi^*(0)$, i.e., the fact that $\bar\psi_n\to \psi^*$ uniformly on every compact set.
 We have thus found a subsequence of solutions $(x^*,\psi_n,\lambda_n)$ of system \rref{sys_x_}--\rref{sys_l_} that converges to  $(x^*,\psi^*,\lambda^*)$ and satisfies \rref{dob_k}.
\bo
\begin{remark}
   In the case of $\lambda^*=0$, condition \rref{620} of {Proposition}~\ref{s3_dd_kab1} can be weakened to
 \begin{equation}\label{641}
 \frac{I(b_n;\tau_n)-I(b^*;\tau_n)}{||I(b^*;\tau_n)||}\to 0 \textrm{ if }
  J(b_n,u^*;\tau_n)-J(b^*,u^*;\tau_n)\to 0.
 \end{equation}
   Indeed, \rref{620} was only necessary for \rref{631}; but in the degenerate case we can take the co-state arc in form \rref{0p}; then, we would have $-\lambda_n I(x_n(0);\tau'_n)=\psi_n(0),||\psi_n(0)||=1$; thus, \rref{631} would follow from \rref{641}.
\end{remark}

\doc {\it of Proposition}~\ref{s3_dd_kab800}.

Set $\zeta=-\frac{\partial l}{\partial x}(b^*).$
We claim that  $I(b^*;\tau_n)\to\zeta$
 as $n\to\infty.$
Proceed by contradiction, then there exist a subsequence $\tau'\subset\tau$
 and $\xi\in{\mathbb{X}}\setminus\{0\}$
with $I(b^*;\tau_n)\xi<\zeta\xi-3$
 for all $n\in{\mathbb{N}}$.
By uniform continuity of $I$ in $G\times{\mathbb{T}}$ and continuity of
$\frac{\partial l}{\partial x}$ at $b^*$,
 there exists a $\bar\delta>0$
  such that
  \[I(b^*+s\xi;\tau_n)\xi<\zeta\xi-2,\ \  b^*+\bar\delta\xi\in G,\ \
  \bar\delta\xi\zeta\leq l(b^*)-l(b^*+\bar\delta\xi)+\delta  \quad\forall s\in(0,\bar\delta], n\in{\mathbb{N}}.\]
By integration, we have
$J(b^*+\bar\delta\xi,u^*;\tau_n)-J(b^*,u^*;\tau_n)<\bar\delta\zeta\xi-2\bar\delta.$
i.e.
\[J(b^*+\bar\delta\xi,u^*;\tau_n)-J(b^*,u^*;\tau_n)+l(b^*+\bar\delta\xi)-l(b^*)<-\bar\delta\]
for all $n\in{\mathbb{N}}$.
By $b^*+\delta\xi\in G\subset{\mathcal{C}}$
 we have contradiction with Definition~\ref{275} of a uniformly
 $\tau$-overtaking optimal process.
Thus, $I(b^*;\tau_n)\to \zeta$
 as $n\to\infty.$

We claim that all assumptions of Proposition \ref{s3_dd_kab1} hold.
Proceed by contradiction, then there exists  a sequence  of $b_n\in{\mathbb{X}}$ such that
\[b_n\to b^*,\quad I(b_n;\tau_n)-I(b^*;\tau_n)\not\to 0,\quad J(b_n,u^*;\tau_n)-J(b^*,u^*;\tau_n)\to 0.\]
By $I(b^*;\tau_n)\to \zeta$,
 passing to the subsequence if necessary, we can say that, for some $\xi\in{\mathbb{X}}\setminus\{0\},$
we have $I(b_n;\tau_n)\xi<\zeta\xi-3, b_n\in G$ for all $n\in{\mathbb{N}}.$
By  continuity of $\frac{\partial l}{\partial x}$
 at $b^*$,
 there exists a natural number  $N$
  such that $l(b_n+\delta\xi)-l(b_n)\leq-\delta\xi\zeta+\delta$
     for $n>N$, $\delta\in(0,1/N].$
By uniform continuity of $\frac{\partial l}{\partial x}$ on $G$,
 there exists a positive number $\bar\delta<1/N$
 such that
  \[I(b_n+s\xi;\tau_n)\xi<\zeta\xi-2,\ \  b_n+\bar\delta\xi\in G,\ \
  \bar\delta\xi\zeta\leq l(b_n)-l(b_n+\bar\delta\xi)+\bar\delta\quad\forall s\in(0,\bar\delta], n\in{\mathbb{N}}.\]
By integration, we have
$J(b_n+\bar\delta\xi,u^*;\tau_n)-J(b_n,u^*;\tau_n)-\bar\delta\zeta\xi<-3\bar\delta,$ 
i.e.
\[J(b_n+\bar\delta\xi,u^*;\tau_n)-J(b_n,u^*;\tau_n)-l(b_n)+l(b_n+\bar\delta\xi)<-2\bar\delta \qquad  \forall n\in{\mathbb{N}}.\]
Then, for sufficiently large $n$, $l(b_n)\to l(b^*)$, $J(b_n,u^*;\tau_n)-J(b^*,u^*;\tau_n)\to 0$
 implies
$J(b_n+\bar\delta\xi,u^*;\tau_n)+l(b_n+\bar\delta\xi)-J(b^*,u^*;\tau_n)-l(b^*)<-\bar\delta.$
By $b_n+\bar\delta\xi\in G\subset{\mathcal{C}}$
we have contradiction with Definition~\ref{275} of a uniformly
$\tau$-overtaking optimal process.
Thus, all assumptions of Proposition~\ref{s3_dd_kab1} are satisfied.

Then, by $N_L^{{\mathcal{C}}}(b^*)=\{0\}$
we have $\lambda>0$ for any $\tau$-limiting solution $(\lambda,\psi)$,
and by \rref{1s} we obtain $\psi(0)=\lambda\frac{\partial l}{\partial x}(b^*)$;
 i.e.,
  accurately to a positive factor, there exists a unique
    $\tau$-limiting solution $(1,\psi^*)$;
    this solution is given by \rref{1ss} with $I_*=-\lambda\frac{\partial l}{\partial x}(b^*).$

In the case of  uniformly overtaking optimal process $(b^*,u^*)$,
by Remark~\ref{11}, for each unbounded sequence of positive numbers $\tau_n$
there exists a
$\tau$-limiting solution,
  i.e.,  accurately to a positive factor, this is $(1,\psi^*)$.
  Since  $I_*$
  independent of~$\tau$,
   we have   \rref{1sss}, for $T=0$, $I_*=-\frac{\partial l}{\partial x}(b^*)$,
   \rref{1sss} implies \rref{900}.
     \bo

 \doc {\it of Proposition}~\ref{s3_dd_kab}

Let us show the statement assuming that
$J(b^*,u^*;\tau_n)=0$ for all $n\in{\mathbb{N}}$.
In this case, put
 $b=b^*$ in \rref{2742}. Now, for all $n,k\in{\mathbb{N}}$ $(n>k)$, we have
\begin{equation}
\label{731}
   ||I(b^*;\tau_n)-I(b^*;\tau_k)||
\leqref{2742}
\omega\left(1/\tau_k,J(b^*,u^*;\tau_n)\!-\!J(b^*,u^*;\tau_k)\right)=\omega(1/\tau_k,0).
\end{equation}
 In particular, there exists a finite limit $\displaystyle I_*=\lim_{n\to\infty}I(b^*;\tau_n).$

 Consider  an arbitrary subsequence of $b_n$ that converges to~$b^*$. By theorem on continuous dependence, we have, as $n\to\infty$,
 $I(b_n;\tau_k)\to I(b^*;\tau_k),$
 $J(b_n,u^*;\tau_k)\to J(b^*,u^*;\tau_k)=0$ for each $k\in{\mathbb{N}}$; now,
\begin{subequations}
\begin{eqnarray}
\lim_{n\to\infty}
 ||I(b^*;\tau_k)-I(b^*;\tau_n)||&\leqref{731}&
 \omega(1/\tau_k,0),\label{732}\\
\lim_{n\to\infty}
||I(b_n;\tau_k)-I(b^*;\tau_k)||&=&0,\label{733}\\
\lim_{n\to\infty}||I(b_n;\tau_n)-I(b_n;\tau_k)||&\leqref{2742}&
\lim_{n\to\infty}\omega\big(1/\tau_k,|J(b_n,u^*;\tau_n)\!-\!J(b_n,u^*;\tau_k)|\big)
\nonumber\\ &=& \lim_{n\to\infty}\omega\big(1/\tau_k,|J(b_n,u^*;\tau_n)|\big)
\label{734}
\end{eqnarray}
\end{subequations}
Summing \rref{732}--\rref{734}, we see, from the triangle inequality, that, for each $k\in{\mathbb{N}}$,
\begin{eqnarray*}
\lim_{n\to\infty}||I(b^*;\tau_n)\!-\!I(b_n;\tau_n)||&\leq& \omega(1/\tau_k,0)+
\lim_{n\to\infty}||I(b_n;\tau_n)\!-\!I(b^*;\tau_k)||\\
&\leqref{2742}&\omega(1/\tau_k,0)+\lim_{n\to\infty}\omega\big(1/\tau_k,|J(b_n,u^*;\tau_n)|\big).
\end{eqnarray*}
 Passing to the limit as $k\to \infty $,  by $J(b^*,u^*;\tau_n)=0$, we have
 \[||I(b_n;\tau_n)\!-\!I(b^*;\tau_n)||\to 0\textrm{ as }J(b_n,u^*;\tau_n)\!-\!J(b^*,u^*;\tau_n)\to 0\]

 Now,  from Proposition~\ref{s3_dd_kab1}, we get the results of Corollary~\ref{affs} for $\psi$. Since the sequence of $I(b^*;\tau_n)$ has the finite limit,
  we obtain \rref{1ss}.

 Let us now prove the general case.
 Consider, along with problem \rref{sys0}--\rref{sysK}, the problem
    \[\textrm{Minimize } l(b)-l(b^*)+\int_{0}^\infty \Big[f_0 (x,u,t) - f_0(x^*(t),u^*(t),t)\Big]dt
 \textrm{ subject to  \rref{sys}--\rref{sysK}}.\]
  These problems clearly share trajectories, optimal processes, system of PMP, and $\tau$-limiting solutions. Condition~\rref{620} also holds for the  new problem, because in this problem, the value of the objective functional is identically zero along $x^*$; thus, in the new problem, formula \rref{1ss} holds for some $\tau$-limiting solution. Therefore it holds in the original problem as well.
 \bo
\begin{remark}
In \cite[Corollary 14]{MyJDC},\cite[Theorem 4]{kab}, existence of a finite limit for $J(b^*,u^*;\tau_n)$ was assumed. As demonstrated by the proof of Proposition~\ref{s3_dd_kab}, this assumption is redundant.
\end{remark}
  \subsection*{Acknowledgements}
   I would like to
express my gratitude to  Ya.V. Salii for
the translation.

\def\theequation{A.\arabic{equation}}
\setcounter{equation}{0}
\section*{Appendix A: a metric for the control system}

  Recall that we can introduce a $\sigma$-metric of space of admissible controls for every $T>0$:
  \[\Delta_T(u,v)\rav \int_0^T 1_{\{s>0\,|\,u(s)\neq v(s)\}}(t)\,dt.\]
  Under this metric, the space of admissible controls becomes a complete metric space. Moreover, the metric can be used to estimate the divergence of trajectories (see, for example, \cite[Lemma~5.1.1]{clarke}).
   The space of admissible controls for infinite interval can also be equipped with such a complete metric (see \cite{norv}), however, it was not clear how to estimate the divergence of trajectories on the whole semi-infinite interval. The variety of metric proposed in this section allows us to do it, at least for a given control system, under assumption of the compactness of set of initial conditions. A variety of such estimate under stronger conditions on a system was proved in \cite{MyArkh} and used in \cite{MyJDC}.

 Let $E$ be a finite-dimensional Euclidean space.
  Consider a map   $a:{\mathbb{U}}\times{E}\times {\mathbb{T}}\mapsto
  E$.
Let us fix a  bounded closed subset set  ${\mathcal{S}}\subset E$ of initial values.
 For each admissible control
$u\in{\mathfrak{U}}$,  consider the differential equation:
\begin{equation}
   \label{a}
   \dot{y}=a(y(t),u(t),t),\qquad
   \forall t\geq 0.
\end{equation}

   We assume that,
   for each admissible control $u\in{\mathfrak{U}}$,
 the map $(y,t)\mapsto a(y,u(t),t)$ is a
 Carath\'eodory map;
 on each bounded subset, a map $(y,u,t)\mapsto a(y,u,t)$
 is integrally bounded and  locally Lipshitz
   continuous on $x$; moreover, each its local solution of \rref{a} can be extended onto the
whole~${\mathbb{T}}$. For every $u\in{{\mathfrak{U}}}$, let us denote the
family of all solutions $y\in C({\mathbb{T}},E)$ of
system~\rref{a} by ${{\mathfrak{A}}}[u]$.

 We also assume that for some  admissible control $u^*\in{\mathfrak{U}}$
 on each bounded subset a map $(y,t)\mapsto a(y,u^*(t),t)$
 is  locally Lipshitz
   continuous on $x$.
Let us fix this admissible control $u^*\in{\mathfrak{U}}$.

 For every point
 $(z,\vartheta)\in{E}\times{\mathbb{T}}$, there exists a unique solution $z\in C_{loc}({\mathbb{T}},E)$
 of the equation
\begin{equation}
   \label{1667}
    \dot{y}=a(y(t),u^*(t),t),\quad y(\vartheta)=z.
\end{equation}
The solution
continuously depends on $(z,\vartheta)$. Let us denote its
initial position~$y(0)$ by $\varkappa(z,\vartheta)$.

 To each $u\in{\mathbb{U}}$, assign the natural number $\lceil u\rceil$, the least natural number not less than $||u||.$
 Note that this function is lower semi-continuous and that it is always $\lceil u\rceil\geq 1.$
 Let us also assign the whole number $\lceil t\rceil\geq 1$ to each time $t\geq 0$.

For all $n\in{\mathbb{N}}$,  consider a bounded subset $G_n$ of $[0,n]\times {\mathbb{X}}$ such that
      \[
      {G}_n\supset\big\{ \left(t,y(t)\right)\,\big|\,y\in{{\mathfrak{A}}}[{u}^*],y(0)\in{\mathcal{S}}, t\in
      {[0,n]}\big\};
      \]
      by the extendability condition for  $a$, this set is bounded.
        Therefore, on this set, the function $a(y,u^*(t),t)$ on ${G}_n$
       is Lipshitz continuous with respect to~$y$
       for the certain Lipshitz constant
       $L_n\in {\mathcal{L}}^1_{loc}([0,n],{\mathbb{T}})$.
       For all $t\geq 0$, define
       \[ M(t)\rav\int_0^t L_{\lceil s\rceil}(s)\,ds.\]

Let us now construct the function $R^a:{\mathbb{N}}\times{\mathbb{T}}\to {\mathbb{T}}.$
For all $t\geq 0,k\in{\mathbb{N}}$,
consider the number
\begin{equation}
\label{nado}
{R}^a(k,t)\rav \sup_{u,v\in {\mathbb{U}},\,\lceil u\rceil,\lceil v\rceil\leq k}\,\sup_{y\in {G}_{\lceil t\rceil}}
 M(t)\big|\big|a(y,u,t)-a(y,v,t)\big|\big|.
\end{equation}
Since~$a$ is integrally bounded on the sets $\{(u,y,t)\,|\,(t,y)\in G_k,||u||\leq k\}$,
 such a number exists for each natural $k$ for almost all $t\in{\mathbb{T}}$.
 Note that ${R}^a(k,t)$ does not decrease as a function of~$k$  for almost all $t\in{\mathbb{T}}$.
 By  \cite[Theorem 3.7]{select}, the map $R^a$ is Borel measurable
 as a function of~$t$  for every $k\in{\mathbb{N}}$.

 Let us now assign the number $w(u,v,t)$ to each pair
  $u,v\in {\mathbb{U}}$ for almost all $t\geq 0$ by the following rule:
\begin{equation*}
 w(u,v,t)
 \rav
 \left\{
 \begin{array} {rcl}        0,       &\mathstrut&     u=v;\\
                            \big\lceil{R}^a(\lceil u \rceil,t)\big\rceil,
                                         &\mathstrut&    u\neq v,\ ||u||\leq||v||;\\
                            \big\lceil{R}^a(\lceil v \rceil,t)\big\rceil,
                                         &\mathstrut&    u\neq v,\ ||u||>||v||.
 \end{array}            \right.
\end{equation*}
 Note that this function is lower-semicontinuous as a function of $u,v$. Indeed, at the point $u=v$ it is true due to ${R}^a\geq 0,$
 and at other points by the monotonicity of ${R}^a$ in $k$ and lower semicontinuity of the mapping $u\mapsto \lceil u \rceil$.
 Then, the function $t\mapsto w(u(t),v(t),t)$ is measurable for each $u,v\in{\mathfrak{U}}$.

  For
all $T\in{\mathbb{T}}$ and $u,v\in{{\mathfrak{U}}}$, let us introduce
  \[\displaystyle \rho(u,v,T)\rav\int_0^T w(u(t),v(t),t)dt\in{\mathbb{T}}.\]
Then, $\rho(u,u,\cdot)\equiv 0$ holds, and, for every
$u\in {{\mathfrak{U}}}$, from $\rho(u,v,T)=0$ for some
$T\in{\mathbb{T}}$, it follows that $u(t)=v(t)$ a.e. on
$[0,T]$.

\begin{lemma}
\label{dop0}
 For every $T>0$, the mapping $(u,v)\mapsto\rho(u,v,T)$ defines a metric on
\[{\mathfrak{U}}_T\rav \{u\in{\mathfrak{U}}\,|\,u(t)=u^*(t) \ \forall t>T\};\]
 under this metric, the space ${\mathfrak{U}}_T$ becomes a complete metric space, and convergence in this metric is no weaker than convergence in measure.

 Moreover, if for some sequence of $u_n\in{\mathfrak{U}}_{\tau_n}$
the sequence of
$\rho(u^*,u_n,\tau_n)$ tends to zero, then the sequence of $u_n$ converges in the metric to $u^*$ on
the whole ${\mathbb{T}}$.
 \end{lemma}

 The next result justifies the introduction of the family $\rho$ of metrics:
\begin{lemma}
\label{dop}
For arbitrary $u\in{{\mathfrak{U}}}$, for every
 solution  $y\in{\mathfrak{A}}[u], y(0)\in{\mathcal{S}}$
 of equation \rref{a}, from $\varkappa(y(t),t)\in{\mathcal{S}}$
 for all $t\in[0, T]$ it follows that
\begin{equation}
\label{1000}
 \big|\big|\varkappa(y(t),t)-y(0)\big|\big|\leq \rho(u^*,u,t) \qquad \forall t\in [0,T].
\end{equation}

 Moreover, for arbitrary $u\in{{\mathfrak{U}}}$, every
 solution  $y\in{\mathfrak{A}}[u], y(0)\in{\mathcal{S}}$ of equation \rref{a} satisfies \rref{1000} if
 $\rho(u^*,u,T)<dist(y(0),bd\,{\mathcal{S}}).$
 \end{lemma}
 \begin{lemma}
\label{dop2}
 For a sequence of $u_n\in{\mathfrak{U}}$ and a
  sequence of  $\td{y}_n\in{\mathfrak{A}}[u_n]$, let
 \[\rho(u^*,u_n,T)\to 0, y_n(0)\to \xi \textrm{ as } n\to\infty \]
 for some $T>0,$ $\xi\in int\,{\mathcal{S}}$.

 Then, the solutions
 $\td{y}_n$ converge to the solution of \rref{a} generated by $u^*$ from the point $\xi$ uniformly in $[0,T]$.
 \end{lemma}
 Thus we construct the estimates (for example, \rref{1000}) of discrepancy between the trajectories generated by arbitrary admissible control $u\in{\mathcal{U}}$ and the trajectories generated by $u^*.$
 \begin{remark}
  The estimates constructed for a pair $(u,u^*)$ may not necessarily hold for an arbitrary pair of controls $(u,u')$; in the general case, to make estimates with respect to a control $u'$ that is different from $u^*$, it may be necessary to construct specific
  $w,\rho;$ however, it is actually possible to provide estimates for the controls from the whole subset of ${\mathcal{U}}$, were  its elements $u'$ to provide the general families of $G_n, L_n$. Moreover, based on the idea from \cite[Sect.~4.3]{MyJDC}, it is also passible to construct a metric that would estimate the discrepancy between trajectories subject to an arbitrary pair of controls $(u,u')\in{\mathcal{U}}\times {\mathcal{U}}.$
   \end{remark}

\subsection{ Proofs of  Lemmata~A.1--A.3.
}
\label{sec:7}
\doc of Lemma~A.1

 Note that since $w(u,v,t)=\max\{{R}^a(\lceil u \rceil,t),{R}^a(\lceil v \rceil,t)\}$ for all $t\in{\mathbb{T}},u,v\in{\mathbb{U}} (u\neq v),$
 $(u,v)\mapsto w(u,v,t)$ is an ultrametric in ${\mathbb{U}}$ for every $t\in{\mathbb{T}}$.
 Hence, $(u,v)\mapsto\rho(u,v,T)$ is a metric of ${\mathfrak{U}}_T$ for every $T\in{\mathbb{T}}$ since it is an integral in $[0,T]$ of a metric.

Moreover, by construction, $w(u,u,t)\equiv 0$, but,
 \begin{eqnarray*}
 w(u,v,t)
 \geq 1\qquad&\forall& \textrm{ a.a. }t>0,u,v\in{\mathbb{U}}(u\neq v);\\
w(u(t),v(t),t)\geq 1_{\{s\geq 0\,|\,u(s)\neq v(s)\}}(t)\qquad &\forall& \textrm{ a.a. }t>0,u,v\in{\mathfrak{U}};\\
\rho(u,v,T)\geq \Delta_T(u,v)=\int_0^T 1_{\{s>0\,|\,u(s)\neq v(s)\}}(t)\,dt\qquad &\forall& T>0,u,v\in{\mathfrak{U}}.
\end{eqnarray*}
 Thus, the introduced metric is not weaker than $\Delta_T$ for every $T>0$. In particular, since $\Delta_T$ turns ${\mathfrak{U}}_T$ into a complete space, $(u,v)\mapsto\rho(u,v,T)$ retains this property. Moreover, since for arbitrary  $u\in{\mathfrak{U}}_T$ and a sequence of $u_n\in{\mathfrak{U}}_T$, the convergence $\Delta_T(u_n,u)\to 0$ is equivalent to the convergence of $u_n|_{[0,T]}$ to $u|_{[0,T]}$ in measure, $\rho(u_n,u,T)\to 0$ also implies that $u_n$ converges in measure to $u$ on the whole ${\mathbb{T}}$.
  \bo

 \doc {\it of Lemma~A.2.}

 Let us now consider
 $u\in{{\mathfrak{U}}},$ $y\in{{\mathfrak{A}}}[u],$
         $T\in{\mathbb{T}}$ such that
          $\varkappa(y(t),t)\in{\mathcal{S}}$ for any $t\in[0,T]$.
           The solution  $z$  of system \rref{a} passes through the point $\left(T,y(T)\right)$. We can assume that $T\in[n-1,n)$ for some $n\in{\mathbb{N}}$.

  For all $t\in[0,T]$, define the number
  $\varrho(t)\rav M(t)||z(t)-y(t)||.$

          By construction of $G_n$, we have
        $\left(t,y(t)\right),\left(t,z(t)\right)\in  G_n$.

        Now, by definitions of $L_n,w$, we have, for almost all
        $t\in[0,T]$,
        \begin{eqnarray*}
           M(t)||\dot{z}(t)\!-\!\dot{y}(t)||&\leq&
            M(t)||a({y}(t),u^*(t),t)\!-\!\dot{y}(t)||
            +M(t)||\dot{z}(t)\!-\!a({y}(t),u^*(t),t)||
             \\
              &\leqref{nado}&
               w(u^*(t),{u}(t),t)+M(t)L_{\lceil t \rceil}(t)||z(t)-y(t)||\\
                &=&
                 w(u^*(t),{u}(t),t)+L_{\lceil t \rceil}(t)\varrho(t).
       \end{eqnarray*}
        Substituting this, we have
        \begin{eqnarray*}
        \frac{d\varrho}{dt}(t)= \frac{d}{dt}\big(M(t)||z(t)-y(t)||\big)&\geq&-
        M(t)||\dot{z}(t)\!-\!\dot{y}(t)||+L_{\lceil t \rceil}(t)\varrho(t)\\
        &\geq&-
         w(u^*(t),{u}(t),t)-L_{\lceil t \rceil}(t)\varrho(t)+L_{\lceil t \rceil}(t)\varrho(t)\\
         &=&-w(u^*(t),{u}(t),t).
       \end{eqnarray*}
       Then, by $\varrho(T)=0$ and the comparison theorems, $\varrho(\cdot)$ does not exceed the solution of the Cauchy problem
       \[ \dot{\varpi}(t)=-w(u^*(t),{u}(t),t) \qquad\varpi(T)=0\]
       in $[0,T]$, i.e., ${\varpi}(t)\rav\rho(u^*,{u},T)-\rho(u^*,{u},t)$.
       Thus,
      $\varrho(t)\leq\rho(u^*,{u},T)-\rho(u^*,{u},t)$ for all $t\in[0,T].$
       In particular, since $z(T)=y(T)$,$\varkappa(z(T),T)=z(0)$,$M(0)=1$,
       we have
       \[||y(0)-\varkappa(y(T),T)||=M(0)||y(0)-z(0)||=\varrho(0)\leq\rho(u^*,{u},T).\]
        Since $t\mapsto\rho(u^*,{u},t)$ is monotonous, the proof of
        \rref{1000} is complete.

 Let us now prove the remaining statement of the Lemma.
 Assume the contrary and let $\rho(u^*,u,T)<dist(y(0),bd\,{\mathcal{S}})$ for some $T>0$, $u\in{{\mathfrak{U}}}$, and
  $y\in{\mathfrak{A}}[u] (y(0)\in{\mathcal{S}})$;  \rref{1000} does not hold.
 Therefore, as proved above, $\varkappa({y}(\vartheta),t)\not\in{\mathcal{S}}$
 holds for some $\vartheta\in[0,T]$ and there exists the greatest
 $T_0\in[0,\vartheta]\subset[0,T]$ such that
 $\varkappa({y}(t),t)\in{\mathcal{S}}$ for all $t\in[0,T_0].$
 As proved above, \rref{1000} holds for all $t\in[0,T_0]$; in particular,
 \[||\varkappa({y}(T_0),T_0)-{y}(0)||\leq \rho(u^*,u,T_0)\leq \rho(u^*,u,T).\]
 By the continuity of $t\mapsto \varkappa({y}(t),t)$, we have
 $\varkappa({y}(T_0),T_0)\in bd\,{\mathcal{S}},$
 i.e.,
 \[\rho(u^*,u,T)<dist(y(0), bd\,{\mathcal{S}})\leq||\varkappa({y}(T_0),T_0)-{y}(0)||,\]
 whence $\rho(u^*,u,T)<\rho(u^*,u,T),$
 which contradicts the assumption. The obtained contradiction completes the proof of Lemma A.2.
 \bo
\doc of {Lemma}~A.3.

 Starting with a certain~$n$, we have
$\rho(u^*,u_n,T)<dist(y_n(0),bd\,{\mathcal{S}})$;
then, estimate \rref{1000} holds for all $t\in[0,T],u=u_n$. Now, since the solution of system \rref{a} generated by the control $u^*$ from the initial position $\varkappa(\td{y}_n(t),t)$ passes through the point $\left(t,\td{y}_n(t)\right)$, the required uniform convergence follows straight from the continuous dependence of solutions of differential equation~\rref{1667} on initial conditions.
    \bo

\section*{Appendix B: proof of Conjecture~\ref{924} for Example 1}
\label{a924}
\def\theequation{B.\arabic{equation}}
\setcounter{equation}{0}

 For every $z\in{\mathbb{R}}$, let $g(z)=z(z^4-5).$ Note the easily verified inequalities
\begin{eqnarray}
  g(z)>z^4\textrm{ \  if \ } z\geq 2;\qquad
  g(z)\geq -5z\textrm{ \ if \ } z\geq0.\label{705}
\end{eqnarray}

 Note that for all $b^*=0,u^*\equiv 0$, in addition to
$x(b^*,u^*;\cdot)\equiv 0$, we have $J(b^*,u^*;\cdot)\equiv 0$.
Then, by definition of uniformly overtaking optimality, it is enough to prove that in this example
\begin{equation}
  \label{72}
 \liminf_{T\to\infty} \inf_{(b,u)\in{\mathcal{C}}\times{\mathcal{U}}} J(b,u;T)\geq 0.
\end{equation}

    Note that for every $T\geq 0$, \rref{72} does reach the infimum on some process $(b_T,u_T)$; denote the trajectory generated by this process by $x_T.$ Without loss of generality, we may assume $u_T|_{(T,\infty)}\equiv 0.$

 Let us show that $b_T\geq 0$ for all $T>0.$
 If $b_T\leq 0$ holds for some $T>0$, then, for an admissible trajectory $y(\cdot)\rav\max\{0,x_T(\cdot)\}$, we have
 $g\big(x_T(t)\big)<g\big(y(t)\big)$ if $t> 0$ and $x_T(t)<0$.
 Then, by optimality of $(b_T,u_T)$ for its own problem, we have $x_T\geq 0$, in particular
 $b_T\geq 0.$

 Let us now show that $b_T<3/2$. Indeed, assume the contrary: let $b_T\geq 3/2.$ Admissible trajectories are nondecreasing, hence $x_T(t)\geq3/2$ for all $t\geq 0$; then,
 $g\left(x_T(t)\right)>0$ for all $t\geq 0$, whence $J(b_T,u_T;T)>0=J(0,0;T)$, which contradicts the optimality of $(b_T,u_T).$

 Thus, $b_T\in[0,3/2]$, i.e., $b_T\in int\, {\mathcal{C}}.$
 Let us find a solution $(\psi_T,\lambda_T)$ of PMP that corresponds to $(b_T,u_T)$.
 Here,
  $\psi_T(0)=0$ by $\lambda_T\upartial_L l(0)+N_L^{{\mathcal{C}}}(0)=\{0\},$
 whence $\lambda_T>0.$
 Moreover, $\psi_T(T)=0.$
 Now, \rref{sys_psi} is equivalent to the equation
 \[\dot{\psi}_T=-\psi_T\frac{\partial f\left(x_T(t)\right)}{\partial x}+\lambda_T \frac{\partial g\left(x_T(t)\right)}{\partial x}=
 -\frac{\psi_T}{2}\max\left(1,x_T(t)\right)+5\lambda_T e^{-2t}(x_T^4(t)-1).\]
 Since $x_T$ is monotonically nondecreasing, $\psi_T$ can change its sign at most twice: once before $x_T(t)$ becomes equal to $1$, and the second time after that. Now, because
  $\psi_T(0)=0=\psi_T(T)$, there indeed exists  a $\vartheta(T)\in(0,T]$ such that
   $x_T\left(\vartheta(T)\right)=1$; moreover
  $\psi_T(t)<0$ in $(0,T),$ whence by \rref{maxH} we have
   $u_T\equiv 0$ in $[0,T].$ Thus, $x_T=x(b_T,0;t)$

  Solving the equation $\dot{z}=f(z)$ with the initial condition $z\left(\vartheta(T)\right)=1,$
  we obtain
  $x_T(s)=\frac{2}{\vartheta(T)+2-s}$ for all $s\in[0,\vartheta(T)].$
  Since $x_T(s)\leq \frac{s+2}{\vartheta(T)+2}$ for  $s\in [0,\vartheta(T)]$,
 we have
   \begin{eqnarray*}
  J\left(b_T,u_T;\vartheta(T)\right)&\geq& \int_{0}^{\vartheta(T)} e^{-2s}g\big(x_T(s)\big)ds\\
  &\geqref{705}&
  -5\int_{0}^{\vartheta(T)} e^{-2s}x_T(s)ds\geq  -\frac{7}{\vartheta(T)+2}.
  \end{eqnarray*}
  Then, it is easy to see that
  $2x_T(s)=e^{s-\vartheta(T)}+1$ for all $s\geq \vartheta(T).$
  Now, for all $s\geq\vartheta(T)+\eta$, we have $g\left(x_T(s)\right)\geq g(\sqrt[4]{5})=0,$
  and  $s\mapsto J(b_T,u_T;s)$ reaches the minimal value at
  $t=
  \vartheta(T)+\eta$, where  $\eta\rav\ln (-1+\sqrt[4]{80})$.
  Then, for all $t>\vartheta(T)$, we have
   \begin{eqnarray*}
  J(b_T,u_T;t)-J\left(b_T,u_T;\vartheta(T)\right)&\geq& J(b_T,u_T;\vartheta(T)+\eta)-J\left(b_T,u_T;\vartheta(T)\right)
  \\
  &\geqref{705}&
  -5\int_{\vartheta(T)}^{\vartheta(T)+\eta} e^{-2s}x_T(s)\,ds\\
  &\geq&-5\int_{\vartheta(T)}^{\vartheta(T)+\eta} e^{-2s} e^{s-\vartheta(T)}\,ds\geq
  -5e^{-2\vartheta(T)}.
  \end{eqnarray*}
  Now, by $T\geq\vartheta(T)$, we have
  \begin{equation}
  \label{73}
  J(b_T,u_T;t)>-\frac{7}{\vartheta(T)+2}-5e^{-2\vartheta(T)}\geq-\frac{12}{\vartheta(T)+2}\qquad\forall t\geq\vartheta(T).
  \end{equation}

 In addition to this, for all $s>\vartheta(T)+2$, we have $x_T(s)\geq x_T\big(\vartheta(T)+2\big)>e>2$. By
  \rref{705},
  $g\big(x_T(s)\big)\geq x^4_T(s)>e^{-2\vartheta(T)+2s},$
    and
\begin{eqnarray}
  J(b_T,u_T;T)&=&
  J\left(b_T,u_T;2+\vartheta(T)\right)+\int_{\vartheta(T)+2}^T e^{-2s}g\big(x_T(s)\big)\,ds\nonumber\\
  &\geq& J\left(b_T,u_T;2+\vartheta(T)\right)+\int_{\vartheta(T)+2}^T e^{-2s} e^{-2\vartheta(T)+2s}\,ds\nonumber\\
   &\geqref{73}& -\frac{12}{\vartheta(T)+2}+(T-\vartheta(T)-2)e^{-2\vartheta(T)}.\label{74}
  \end{eqnarray}

  Let us finally prove \rref{72}.
  Assume the contrary; then, there exist a number $\delta>0$ and an unbounded increasing sequence
  of $T(n)\in{\mathbb{T}}$ such that $J\left(b_{T(n)},u_{T(n)};T(n)\right)<-\delta$ for all $n\in{\mathbb{N}}$.
  Passing to the subsequence, we can rest assured that either the sequence $\vartheta\left(T(n)\right)$ is unbounded and increasing or it is bounded.
  In the first case, the contradiction is obtained by passing in \rref{73} to the limit as $T(n)\to\infty$, while in the second, the contradiction surfaces in \rref{74}.
  Thus, \rref{72} is proved, and Conjecture~\ref{924} with it.
\end{document}